\documentclass[a4paper,10.5pt]{amsart}

\usepackage{amsmath}
\usepackage{amssymb}
\usepackage{amsthm}
\usepackage{amscd}
\usepackage[dvips]{graphicx}
\usepackage{enumerate}
\usepackage{breqn}
\usepackage{cleveref}
\usepackage{enumerate}
\usepackage{xcolor}

\newtheorem{proposition}{Proposition}[section]
\newtheorem{theorem}{Theorem}[section]
\newtheorem{lemma}{Lemma}[section]

\numberwithin{equation}{section}
\newcommand{\real}{\mathbb{R}}

\def\eqn {\begin{equation}}
\def\eeqn {\end{equation}}

\def\real{{\mathbb R}}

\def\ep{\epsilon}

\def\pa{\partial}

\def\A{\mathcal A}

\def\F{\mathcal F}
\def\K{\mathcal K}
\def\L{\mathcal L}

\def\N{\mathcal N}
\def\O{\mathcal O}

\def\S{\mathcal S}

\begin{document}
\title{Global Continuation of a Vlasov Model of Rotating Galaxies}
\author{Walter A. Strauss} 
\address{Department of Mathematics, Brown University, Providence, RI 02912}
\email{walter_strauss@brown.edu}
\author{Yilun Wu}
\address{Department of Mathematics, University of Oklahoma, Norman, OK 73019}
\email{allenwu@ou.edu}
\date{}


\begin{abstract}
A typical galaxy consists of a huge number of stars attracted to each other by gravity.  For instance, the Milky Way has about $10^{11}$ stars.  
Thus it is typically modeled by the Vlasov-Poisson system.  
We prove an existence theorem for axisymmetric steady states of galaxies that may rotate  rapidly.  Such states are given in terms of a fairly general function $\phi$ of the particle energy and angular momentum.  
The set $\K$ of such states form a connected set in an appropriate function space.
Along the set $\K$, we prove under some conditions that either (a) the supports of the galaxies become unbounded 
or (b) both the rotation speeds and the densities somewhere within the galaxy become unbounded.  
\end{abstract}

\maketitle

\section{Introduction}

We consider a continuum of particles with a fixed total mass that are 
attracted to each other by gravity but subject to no other forces.  
Initially they are static and spherical, but if they begin to rotate 
around a fixed axis, they flatten at the poles and expand at the equator.  
This is a very simple way to think about a rotating collection of particles.  
In this paper we model this configuration of particles 
in the standard way by the Vlasov-Poisson system (VP).  
It is the standard model of stellar systems such as galaxies \cite{binney2011galactic}. 
We look for for steady states of the resulting configuration  
and find a connected set $\K$ of such states with constant mass 
and possibly large rotation speeds.

\bigskip 
The VP system is a continuous model for the particle distribution $f(x,v)$ 
in phase space $(x,v)$ where $x\in\real^3$ is position and $v\in\real^3$ is velocity.  
The Vlasov equation (also called the collisionless Boltzmann equation) is 
\eqn\label{eq: Vlasov}
\pa_tf+ v\cdot\nabla_x f+\nabla_xU\cdot \nabla _v f=0,       \eeqn
where $f(x,v)\ge0$ and $\nabla_x U$ is the force of gravity.  
The gravitational potential $U$ is 
\eqn\label{eq: Poisson}
U=\frac{1}{|\cdot|}*\rho,\quad \rho(x)=\int_{\real^3}f(x,v)~dv, 
\eeqn
where $\rho$ is the macroscopic particle density in physical space. 
In general, the macroscopic velocity $V(x)=\int_{\real^3}v\ f(x,v)\ dv$ does not vanish.  

In this paper we only consider solutions that are independent of time $t$ 
and are rotating around the $x_3$ axis.   
Jeans' Theorem \cite{binney2011galactic} states that every steady state depends only on the 
invariants of the flow.  Two of these invariants are the energy and the 
$x_3$ component of the angular momentum.  
We only consider solutions that are axisymmetric and of the form 
\eqn         \label{eq: f form}
f(x,v)=\phi\left(\tfrac12|v|^2-U(x)-\alpha,~\kappa(x_1v_2 - x_2v_1)\right),       \eeqn
where 
$\kappa$ is a measure of the intensity of rotation and 
the constant $\alpha$ is included for convenience.  
The microscopic density function $\phi$ is nonnegative.  
Then $f$ is itself an invariant of the flow and so it is automatically a solution of the 
Vlasov equation \eqref{eq: Vlasov}.

Let $C_c(\real^3)$ be the Fr\'echet space of axisymmetric continuous functions 
on $\real^3$ with compact support.  
Its topology, defined by seminorms, is that of uniform convergence on each compact set.  
By a {\it solution} of our problem we mean a triple 
$(\rho,\alpha,\kappa) \in C_c(\real^3)\times\real\times\real$ 
that satisfies \eqref{eq: Poisson} and \eqref{eq: f form}.  

Our main result is as follows.  
\begin{theorem}      \label{main thm}
Let $(\rho_0,\alpha_0)$ be a non-rotating $(\kappa=0)$ spherically symmetric solution 
be given (as described in Lemma \ref{lem: non-rot} below).  
Let $M_0=\int_{\real^3} \rho_0 dx$ be its total mass.    
Assume that a microscopic density function $\phi(E,L)$ is given that 
satisfies the assumptions \eqref{phi: pos} through \eqref{phi: upper bound}.  
Then there exists a set $\K$ of solutions $(\rho,\alpha,\kappa)$ with $\rho\ge0$ 
and $\rho\in C^1$ that satisfies the following properties. 

$\bullet \ \ \K$ is a connected set in the space $C_c(\real^3)\times\real\times\real$.  

$\bullet\ \ \K$ contains the non-rotating solution.

$\bullet$\ \ All the elements in $\K$ have the same total mass $M_0$.  

$\bullet$\ \ Either 
\eqn   \label{support unbdd}
   \sup\left\{ |x|\ \Big| \ \rho(x)>0, \ (\rho,\alpha,\kappa)\in\K\right\}=\infty   \eeqn
or 
\eqn     \label{density unbdd}
\sup\left\{ \rho(x) \ \Big| \ x\in\real^3, (\rho,\alpha,\kappa)\in\K\right\}=\infty    .  
\eeqn
\end{theorem}

\begin{theorem}      \label{main thm2} 
Under the same setup as in Theorem \ref{main thm}, if $\phi(E,L)$ satisfies assumptions \eqref{phi: pos} through \eqref{cond: increasing mass} and \eqref{phi: upper bound} but with the more restrictive range of $\Lambda$: $0<\Lambda<4$, then the same conclusion as in Theorem \ref{main thm} is valid 
provided \eqref{density unbdd} is replaced by 
\eqn     \label{speed unbdd} 
 \sup\left\{|\kappa| \ \Big|\  (\rho,\alpha,\kappa)\in\K\right\}=\infty .      \eeqn
\end{theorem}

Alternative \eqref{support unbdd}  means that $\rho$ has unbounded support.  
Because the total mass is fixed, it would mean that there can be a very large 
set in $\real^3$ on which $\rho$ is very small.  
Alternative \eqref{density unbdd} means that the macroscopic density is unbounded.  
The other alternative \eqref{speed unbdd} means that the rotation speed is unbounded.  
Which of these alternatives can actually happen is an open question.  

As discussed later, a {\it simple example} of $\phi$ that satisfies the conditions 
of Theorem \ref{main thm} is $\phi(E,L) = (E)_-^\nu~p(L)$, 
where $(E)_-$ denotes the negative part of $E$, $\nu\in(-\frac12,\frac72)$ , $\nu\ne\frac32$, 
and $p$ is any positive polynomial.  The same function  satisfies the conditions of 
Theorem \ref{main thm2} provided $p$ is a positive quadratic polynomial.   
More examples are given at the end of Section 3.  

\bigskip 
There is a voluminous classical literature on fluid models of rotating stars 
going back notably to MacLaurin, Poincar\'e and Lichtenstein.  
On the other hand, the Vlasov system with gravity is a more natural model of galaxies.  
Spherical solutions of this system are well-known in the mathematical literature 
(see \cite{rein2007collisionless}) but non-spherical ones have not been studied much.   
See \cite{rein2007collisionless} for a nice survey of the initial-value problem.  

The first mathematical treatment of steady states of axisymmetric solutions 
of VP that we are aware of is by Rein \cite{rein2000compact}, 
who constructs a local curve of nearly spherical solutions 
by means of the implicit function theorem.  
A similar theorem may be found in \cite{jang2017slowly}.  
Such a local curve is also constructed by Andreasson {\it et al.} 
within the context of general relativity in \cite{andreasson2014rotating}.  
In \cite{strauss2017steady}, Sections 6-8, we have constructed such a local curve 
with the additional property that the mass remains constant 
(independent of the rotation speed).  
On the other hand, Guo and Rein \cite{rein2003stable} construct stable steady states by minimization of an energy-Casimir functional for a certain class of $\phi$.  
The stability of spherical states has also been studied by Guo and Lin \cite{guo2008unstable}.  

The point of the present paper is to construct a ``global" connected set of steady states 
that have compact support and constant mass 
and that emanate from a spherical solution but deviate far from it.  
The constant mass condition means that there is no loss or gain of particles 
as the star changes its rotation speed.   
Our method has much in common with \cite{strauss2019rapidly}, where a global curve 
is constructed for a fluid model of a rotating star.  However, the two different models 
require quite different analyses in some respects.  

\bigskip
In the next section we state our detailed assumptions.  
The remainder of the paper is devoted to the proof of the two theorems.  
In Section 3 we provide details about the non-rotating spherical solution and its mass.  
Section 4 provides the proof of global continuation by means of a global implicit function theorem.  
Sections 5 and 6 are devoted to the proof that either the support of $\rho$ is unbounded 
or $\rho$ itself is unbounded (Theorem \ref{main thm}).  
This is the most subtle part of our paper.  
In order to elucidate the proof, it is first proven for a special case 
and only later for a general case.  Key points in the proof are the use of 
a Gagliardo-Nirenberg inequality and of the constant mass condition.  
Section 7 is devoted to proving that either the support of $\rho$ is unbounded or the rotation speed is unbounded (Theorem \ref{main thm2}).    

\section{Setup and Basic Assumptions}\label{sec: setup}
It is easily verified that for any $(\kappa,\alpha)\in\real^2$, $U\in C^{2}(\real^3)$ axisymmetric  and $\phi\in L^1_{\text{loc}}(\real^2)$, 
$f$ given by \eqref{eq: f form} 
is a weak solution to \eqref{eq: Vlasov} in the sense that 
\eqn
f(x,v)=f(\psi(x,v,t))  ,  \eeqn
where $\psi(x,v,t)$ is the flow map of the vector field $(v,\nabla U(x))$ on $\real^6$. 
Given a function $\phi$ on $\real^2$, we {\it define}  
\eqn
w(\kappa,r,u)=\int_{\real^3}\phi\left(\frac12|v|^2-u,\kappa(x_1v_2-x_2v_1)\right) \, dv,\quad r=(x_1^2+x_2^2)^{1/2}.    \eeqn
Then the Poisson equation \eqref{eq: Poisson} takes the form 
\eqn     \label{eq: Poisson2}
\rho = w\left(\kappa,r,\frac1{|\cdot|}*\rho+\alpha\right).      \eeqn
Given $\phi$, it suffices to solve \eqref{eq: Poisson2} for $\rho, \alpha$ and $\kappa$.  

In order to construct a global set of solutions with fixed total mass, we therefore look for the zero set of the mapping $\F(\rho,\alpha,\kappa)=(\F_1,\F_2)$, where
\begin{equation}\label{eq: F1}
\mathcal F_1(\rho,\alpha,\kappa)= \rho(x) - w\left(\kappa, r, \frac1{|\cdot|}*\rho(x) + \alpha\right),
\end{equation}
\begin{equation}\label{eq: F2}
\mathcal F_2(\rho,\alpha,\kappa) = \int_{\mathbb R^3} \rho(x)~dx- M_0.
\end{equation}
Here $M_0$ is the fixed total mass of the solution set, which we will discuss in more detail in Section 3.   {\it We will always assume that $\rho$ is axisymmetric and even in $x_3$.}

Now we specify the 
{\it assumptions on the prescribed microscopic density function} 
$\phi:\real^2\to \real$. We divide them into several groups:
\begin{enumerate}[(I)]

\item 
We begin with its positivity and regularity properties, 
which will be used to show the regularity and positivity of $\rho$:
\eqn\label{phi: pos}
\phi(E,L)>0 \text{ when }E<0,\quad \phi(E,L)=0 \text{ when }E>0.
\eeqn
\eqn\label{phi: reg}
\phi\big|_{E<0}\in C^1_{\text{loc}}\left((-\infty,0)\times \real\right).
\eeqn

\item 
Next are the properties that will be used to obtain the non-rotating ($\kappa=0$) solutions:
\eqn\label{cond: phi(E,0) 2}
\lim_{E\to-\infty} |E|^{1/2}\phi(E,0)=\infty.    \eeqn
\eqn\label{cond: phi(E,0) 3}
\lim_{E\to-\infty} |E|^{-7/2}\phi(E,0)=0.    \eeqn

\item 
We also require the following assumption, stated rather informally:
\begin{align}\label{cond: increasing mass}
&\text{The total masses of the non-rotating solutions depend} \\
&\text{strictly monotonically on the center densities }\rho(0).\notag
\end{align}
 This property will be used to generate a local curve of solutions emanating from a non-rotating one.   A more precise statement of it is given in \eqref{cond: M 1}, \eqref{cond: M 2}.  

\item The final pair of properties are a lower bound and an upper bound that will be used mainly to narrow down the types of blow up behavior that could occur in the global solution set in Theorems \ref{main thm} and \ref{main thm2}:
\eqn\label{phi: lower bound}
\liminf_{E\to 0^-,|L|\to\infty}|E|^{-\delta}|L|^{-\Gamma}\phi(E,L)>0 
\eeqn
for some $\delta>0$ and $ \Gamma>0$.
Furthermore, for every $B>0$, there exist $C>0$, $\Lambda>0$, $\mu<\tfrac12$ such that 
\eqn\label{phi: upper bound}
\phi(E,L)+|\partial_L\phi(E,L)|\le C|E|^{-\mu}(1+|L|)^\Lambda 
\eeqn
for all $-B<E<0$ and all $L$.


\end{enumerate}
In Proposition \ref{thm: mass cond}, we will provide some simple conditions 
on the density function $\phi$ under which \eqref{cond: increasing mass} is satisfied.  

\section{Non-rotating Solutions}\label{sec: non rot}
We now discuss the non-rotating solutions ($\kappa=0$) in more detail. Using spherical coordinate substitution with polar axis pointing in the direction of $(-x_2,x_1,0)$, and reparametrizing using $E=\tfrac12v^2-u$, $s=\frac{x_1v_2-x_2v_1}{r}$, we can write 
\begin{align}
w(\kappa,r,u)&=\int_{\real^3}\phi\left(\frac12v^2-u,\kappa(x_1v_2-x_2v_1)\right)~dv\notag\\
&=2\pi \int_{-u}^0 \int_{-\sqrt{2(E+u)}}^{\sqrt{2(E+u)}} \phi(E,\kappa r s) ~ds~dE.\label{eq: w}
\end{align}
We easily see from \eqref{eq: w} that $w(\kappa,r,u) = w(-\kappa, r, u)$. The switch from $\kappa$ to $-\kappa$ corresponds to reversing the direction of macroscopic velocity in the fluid. This symmetry means the model is indifferent to the direction of rotation, as in the Euler-Poisson model studied in \cite{strauss2019rapidly}.
We have the following basic properties of $w$.  
\begin{lemma}\label{lem: w basic}
$w\in C^{1}_{\text{loc}}(\real\times [0,\infty)\times \real)$. $w(\kappa,r,u)>0$ if $u>0$, and $w(\kappa,r,u) = 0$ if $u\le 0$. If we define 
\eqn
G(u) = w(0,0,u) = w(0, r, u), 
\eeqn
then $G'(u)>0$ for $u>0$ and $G$ has the limits 
\begin{equation}\label{G: limits}
\lim_{u\searrow 0} u^{-1}G(u) = 0, ~\lim_{u\to\infty} u^{-1}G(u) = \infty,~\lim_{u\to\infty}u^{-5}G(u) = 0.
\end{equation}
\end{lemma}
\begin{proof}The proof resembles that of Lemma 6.4 and Lemma 8.1 of \cite{strauss2017steady}, but we include it here because the hypotheses are a bit different.  
The positivity property of $w$ follows immediately from \eqref{phi: pos}.
By \eqref{phi: upper bound} and the dominated convergence theorem, we differentiate under the integral sign to get
\eqn
\partial_\kappa w(\kappa,r,u)=2\pi r\int_{-u}^0 \int_{-\sqrt{2(E+u)}}^{\sqrt{2(E+u)}}\partial_L\phi(E,\kappa rs)s~ds~dE.
\eeqn

\eqn
\partial_r w(\kappa,r,u)=2\pi \kappa \int_{-u}^0 \int_{-\sqrt{2(E+u)}}^{\sqrt{2(E+u)}}\partial_L\phi(E,\kappa rs)s~ds~dE.
\eeqn
The continuity of both of these functions follows from the dominated convergence theorem. Assuming for the moment that $u>0$, we also have 
\eqn\label{eq: D_u w}
\partial_u w(\kappa,r,u)=\pi\sqrt{2}\int_{-u}^0\frac{\phi\left(E,\kappa r\sqrt{2(E+u)}\right)+\phi\left(E,-\kappa r \sqrt{2(E+u)}\right)}{\sqrt{E+u}}~dE.
\eeqn
It is not hard to see that $\partial_u w(\kappa,r,u)$ is continuous and positive for $u>0$.
Moreover, by \eqref{phi: upper bound} we get for $0<u<1$ that
\eqn
|\partial_u w(\kappa,r,u)|\le C\int_{-u}^0\frac{|E|^{-\mu}}{\sqrt{E+u}}~dE\le Cu^{\frac12-\mu}.
\eeqn
From this inequality and the fact that $w(\kappa,r,u)=0$ when $u\le 0$, we see that $\partial_uw(\kappa,r,0)$ exists and is zero. Also, $\partial_u w$ is continuous everywhere.

We now turn to the study of 
\eqn\label{eq: G}
G(u):=w(0,0,u)=w(0,r,u)=4\pi\sqrt2\int_{-u}^0\phi(E,0)\sqrt{E+u}~dE.
\eeqn
The first limit in \eqref{G: limits} is just a restatement of the fact that $G'(0)=0$. To obtain the second limit, we estimate
\begin{align*}
u^{-1}G(u) &\ge Cu^{-1}\int_{-u}^{-u/2}\phi(E,0)\sqrt{E+u}~dE\\
&\ge C\left(\inf_{E<-u/2}|E|^{1/2}\phi(E,0)\right)\cdot u^{-1}\int_{-u}^{-u/2}|E|^{-1/2}\sqrt{E+u}~dE\\
&= C\left(\inf_{E<-u/2}|E|^{1/2}\phi(E,0)\right)\cdot \int_{-1}^{-1/2}|E|^{-1/2}\sqrt{E+1}~dE
\end{align*}
Thus the second limit in \eqref{G: limits} follows from \eqref{cond: phi(E,0) 2}. To get the last limit, we estimate for any $a>0$:
\begin{align*}
0\le u^{-5}G(u)&\le u^{-5}\int_{-a}^0 \phi(E,0)\sqrt{E+u}~dE+u^{-5}\int_{-u}^{-a}\phi(E,0)\sqrt{E+u}~dE\\
&\le u^{-5}\sqrt u \int_{-a}^0\phi(E,0)~dE+\left(\sup_{E<-a}|E|^{-7/2}\phi(E,0)\right)\cdot u^{-5}\int_{-u}^{0}| E|^{7/2}\sqrt{E+u}~dE\\
&\le \frac{\int_{-a}^0\phi(E,0)~dE}{u^{9/2}}+C\sup_{E<-a}|E|^{-7/2}\phi(E,0).
\end{align*}
It follows that 
\eqn
0\le \limsup_{u\to\infty} u^{-5}G(u)\le C\sup_{E<-a}|E|^{-7/2}\phi(E,0).
\eeqn
The last limit in \eqref{G: limits} follows by sending $a$ to infinity and using \eqref{cond: phi(E,0) 3}.
\end{proof}

We summarize the basic properties of non-rotating solutions as follows.  
\begin{lemma}\label{lem: non-rot} 
For every $R>0$, there exists a radial function $\rho\in C^1(\real^3)$ with $\rho>0$ on $B_R$, $\rho=0$ on $\real^3\setminus B_R$, and a constant  $\alpha<0$ such that $\F_1(\rho,\alpha,0)=0$ on $\real^3$. Furthermore, every nonzero compactly supported continuous solution to $\F_1(\rho,\alpha,0)=0$ satisfies the same conditions.
\end{lemma}
\begin{proof}
Consider the equation
\eqn\label{eq: radial elliptic}
\Delta u +4\pi G(u)=0
\eeqn
in $B_R$ with $u=0$ on $\partial B_R$.
By \cite{ambrosetti1973dual} and \cite{de1982priori}, \eqref{eq: radial elliptic} has a positive $C^2$ solution on $B_R\subset \real^3$ with zero boundary data if $G$ maps $\real^+$ to $\real^+$, is $C^1$ on $\overline{\real^+}$, $G(0)=0$, and satisfies the last two limits in \eqref{G: limits}. By Theorem 1 in \cite{gidas1979symmetry}, $u$ is radially symmetric and strictly decreasing:
\eqn\label{eq: radial decrease}
x\cdot \nabla u(x)<0~\text{ for }0<|x|\le R.
\eeqn
Now define $\rho=G(u)$ on $B_R$ and $\rho=0$ on $\real^3\setminus B_R$. 
The first limit in \eqref{G: limits} implies $\rho\in C^1(\real^3)$. 
By \eqref{eq: radial elliptic}, $u-\tfrac{1}{|\cdot|}*\rho$ is harmonic in $B_R$. 
Since it has a constant boundary value $\alpha=-\tfrac{1}{|\cdot|}*\rho$ on $\partial B_R$, 
it must be constant in all of $B_R$. In other words, we have
\eqn\label{eq: F=0 radial}
\rho-G\left(\frac1{|\cdot|}*\rho+\alpha\right)=0
\eeqn
in $\overline{B_R}$. \eqref{eq: F=0 radial} is $\F_1(\rho,\alpha,0)=0$ by definition. 
To see that \eqref{eq: F=0 radial} holds in $\real^3\setminus B_R$ as well, we only need to show $\tfrac{1}{|\cdot|}*\rho+\alpha<0$ in $\real^3\setminus B_R$, 
since $G(s)=0$ for $s\le 0$. This  follows because $\tfrac{1}{|\cdot|}*\rho+\alpha$ is harmonic outside $ \overline{B_R}$, zero on $\partial B_R$ and negative at infinity.

Now we assume a non-trivial $\rho$ is compactly supported and continuous, with axisymmetry and even symmetry in $x_3$, and $\alpha$ is a constant such that $\F_1(\rho,\alpha,0)=0$. We have $\rho\ge 0$, $\rho\in C^1(\real^3)$ as $G$ is non-negative and $C^1$. 
Let $v=\tfrac1{|\cdot|}*\rho$.   Then $v>0$ and  $\Delta v=-4\pi\rho=-4\pi G(v+\alpha)$.   
Moreover $\alpha<0$, for otherwise $\rho$ could not be compactly supported. 
Applying Theorem 4 and Proposition 1 in \cite{gidas1979symmetry}, 
we conclude that $v$ is radially symmetric about some point in $\real^3$, 
and the same is true for $\rho$. 
We assumed that $\rho$ is axisymmetric and even in $x_3$. It is not hard to see that if $\rho$ is also radially symmetric about a point different from the origin, then it could not be compactly supported unless it is identically zero. Let $\overline{B_R}$ denote the support of $\rho$.  
Letting $u=v+\alpha$, we have $\rho=G(u)$ and $\Delta u \le 0$ in $B_R$, 
and $u=0$ on $\partial B_R$. By the strong maximum principle, $u> 0$ on $B_R$ 
and the same is true for $\rho$. Thus $\rho$ is positive in $B_R$ and zero outside, 
and so it is one of the solutions constructed above.
\end{proof}

Having classified all the non-rotating solutions, we can parametrize them by their center densities.
\begin{lemma}\label{lem: radial uniqueness}
Denote by $\mathcal A$ the set of center values $\rho(0)$ of all the compactly supported non-rotating solutions given in Lemma \ref{lem: non-rot}.   Then 
$\mathcal A$ is an open subset of $(0,\infty)$. For every $a\in \mathcal A$, there exist a unique $\rho(x;a)\in C_c(\real^3)$ and a unique $\alpha(a)$ satisfying $\F_1(\rho(x;a),\alpha(a),0)=0$, $\rho(0;a)=a$. \end{lemma}
\begin{proof}
As in the proof of Lemma \ref{lem: non-rot}, we let $u=\tfrac1{|\cdot|}*\rho+\alpha$, so that  
$\rho=G(u)$. Since $G$ maps $\real^+$ to $\real^+$ diffeomorphically, we can equivalently parametrize the solution by $u(0)$. By radial symmetry, $u$ satisfies the ODE
\eqn\label{eq: radial ODE}
u''+\frac2{|x|}u'+4\pi G(u)=0,
\eeqn
where $'$ denotes the radial derivative. 
We easily see that it has a unique solution $u(r;a)$ satisfying $u(0;a)=G^{-1}(a)$, $u'(0;a)=0$, for any $a>0$. Such a solution must do either of the following:
\begin{enumerate}[(i)]
\item There exists an $R(a)>0$ such that $u(x;a)>0$ for $|x|<R(a)$, and $u(R(a);a)=0$.
\item $u(x;a)>0$ for all $x$.
\end{enumerate}
Since $\mathcal A$ is the set of center values of compactly supported solutions
$\rho$, it follows that $a\in \mathcal A$ if and only if case (i) occurs for $u(x;a)$.  
For any such $a$, by \eqref{eq: radial decrease}, we must have $u'(R(a);a)<0$.  
By the implicit function theorem, for every $b$ in a small neighborhood of $a$, 
there exists $R(b)$ close to $R(a)$ such that $u(R(b);b)=0$. 
This shows that $\mathcal A$ is open.
\end{proof}

Denote the total mass as 
\eqn
M(a)=\int_{\real^3}\rho(x;a)~dx.
\eeqn
It is easy to see that $M\in C^1(\mathcal A)$.  
Choosing any $a_0\in\mathcal A$, we denote $\rho_0(x)=\rho(x;a_0)$, $\alpha_0=\alpha(a_0)$, $M_0=M(a_0)$. In the following discussion we will construct a solution set emanating from $(\rho,\alpha,\kappa)=(\rho_0,\alpha_0,0)$. To that end we replace \eqref{cond: increasing mass} by the following more precise conditions on the mass function $M(a)$:
\eqn\label{cond: M 1}
M'(a_0)\ne 0,
\eeqn
\eqn\label{cond: M 2}
M(a)\ne M(a_0) \text{ for all }a\in \mathcal A, ~a\ne a_0.
\eeqn
Note that \eqref{cond: M 1} and \eqref{cond: M 2} are actually a little weaker than \eqref{cond: increasing mass}. We essentially only need the total mass of the other non-rotating solutions to differ from that of $\rho_0$. 
For instance, if we are in the simple scenario where $\mathcal A$ consists of a single interval, and $M'\ne 0$ on $\mathcal A$, then \eqref{cond: M 1} and \eqref{cond: M 2} are satisfied for any choice of $a_0$.
As was alluded to earlier, the mass conditions \eqref{cond: M 1} and \eqref{cond: M 2} are not directly verifiable on the structure function $\phi$. We provide the following theorem supplying some sufficient conditions on $\phi$ under which \eqref{cond: M 1} and \eqref{cond: M 2} are satisfied.  


\begin{proposition}\label{thm: mass cond}
We have both $\mathcal A=(0,\infty)$ and $M'\ne 0$ on $\mathcal A$  in both of the following cases.  
\begin{enumerate}[(a)]
\item $\phi(E,0)=(E)_-^{\nu}$ for some $\nu\in(-\tfrac12,\tfrac72)$, $\nu\ne\tfrac32 $.
\item $-\tfrac12\phi(E,0)<E\partial_E\phi(E,0)\le \tfrac12\phi(E,0)$ for all $E<0$, 
in addition to the conditions stated in Section \ref{sec: setup}. 
\end{enumerate}
\end{proposition}

\begin{proof}
We first consider the case $(a)$ when $\phi(E,0)=(E)_-^\nu$, then $G(u)=Cu_+^{3/2+\nu}$ for some constant $C>0$. In this case all the non-rotating solutions are related to each other by scaling, as is shown by the explicit scaling symmetry of \eqref{eq: radial elliptic}.  
Therefore $\mathcal A=\real^+$ and 
\eqn\label{M: power}
M(a)=C_0 a^{\frac{3-2\nu}{4+4\nu}}
\eeqn
for some constant $C_0>0$. The range $-\tfrac12<\nu<\tfrac72$ is required 
for $G$ to be $C^1$ and for the radial solutions to be compactly supported. 
The exclusion of $\nu=\tfrac32$ is required in order that $M'\ne 0$. 

In the case of $(b)$, we claim that $G(u)$, given by \eqref{eq: G},   satisfies 
\eqn\label{G: mass cond}
G(u)<uG'(u)\le 2G(u).
\eeqn
In fact, writing $\partial_E(2E\sqrt{E+u})=3\sqrt{E+u}-\frac{u}{\sqrt{E+u}}$ 
and integrating by parts, we get 
\begin{align}
G'(u)&=2\sqrt 2 \pi\int_{-u}^0 \frac{\phi(E,0)}{\sqrt{E+u}}~dE\notag\\
&=\frac{2\sqrt2\pi}u\int_{-u}^0\left(3\phi(E,0)+2E\partial_E\phi(E,0)\right)\sqrt{E+u}~dE.\label{eq: G'}
\end{align}
Comparing \eqref{eq: G'} with \eqref{eq: G}, we see that the inequality in case (b) 
implies \eqref{G: mass cond}. Using \eqref{G: mass cond}, we can employ exactly 
the same argument as in Lemma 2.1 of \cite{strauss2019rapidly} (see also Lemma 4.9 in \cite{strauss2017steady}) with $h^{-1}$ replaced by $G$ to prove that $\mathcal A=\real^+$  and that $M'>0$ on $\mathcal A$. 
\end{proof}

{\it Examples of $\phi$:}  \begin{enumerate}[($i$)]   

\item  $\phi(E,L) = (E)_-^\nu ~p(L)$, 
where $(E)_-$ denotes the negative part of $E$, $\nu\in(-\frac12,\frac72)$ , $\nu\ne\frac32$ 
and $p$ is any positive polynomial.  This follows directly from Proposition \ref{thm: mass cond}(a).  

\item  $\phi(E,L) = [(E)_-^{\nu_1} + (E)_-^{\nu_2}]~p(L)$, 
where $\nu_1,\nu_2\in(-\frac12, \frac12)$ and $p$ is any positive polynomial.  
This follows directly from Proposition \ref{thm: mass cond}(b). 

\item  $\phi(E,L) = \left(A + \frac{\sin E}E \right)~p(L)$ for $E<0$, $\phi=0$ for $E>0$, 
$A$ is a large constant and $p$ is any positive polynomial.  
This also follows easily from Proposition \ref{thm: mass cond}(b). 

\item For Theorem \ref{main thm2}
we take the polynomial $p$ to be quadratic.  

\end{enumerate}


\section{Local and Global Continuation}
We employ a convenient function space on which $\F$ is well defined. Let
\eqn
X = \{f:\real^3\to\real~|~f \in C(\real^3) \text{ is axisymmetric, even in }x_3,~\|f\|_X<\infty\},
\eeqn
where 
\eqn\label{eq: X norm}
\|f\|_X=\sup_{x\in\real^3}\langle x\rangle ^4|f(x)|.
\eeqn
Here $\langle x\rangle = (1+|x|^2)^{1/2}$. The power $4$ in \eqref{eq: X norm} could be replaced by any real number bigger than 3 without making any essential changes to the arguments below.
We will show in the following discussion that $\F=(\F_1,\F_2)$ defined by \eqref{eq: F1}, \eqref{eq: F2} is a $C^1$ mapping from $X\times \real\times \real$ to $X \times \real$ and study its zero set.
For $N\in \mathbb N$, let 
\begin{equation}
\mathcal O_N = \left\{(\rho,\alpha,\kappa)\in X \times \mathbb R^2~|~ \alpha<-\frac1N\right\}.
\end{equation}
\eqn
\O_\infty=\bigcup_{N=1}^\infty\O_N=\{(\rho,\alpha,\kappa)\in X\times\real^2~|~\alpha<0\}
\eeqn
By the fact that $w(\kappa,r,u) = 0$ for $u\le 0$ in Lemma \ref{lem: w basic}, we have a simple but important compact support lemma.   
\begin{lemma}\label{lem: comp supp}
For all $(\rho,\alpha,\kappa)\in \mathcal O_N$, $w\left(\kappa, r, \frac1{|\cdot|}*\rho(x) + \alpha\right)\in C(\real^3)$ and is supported in the ball $|x|\le C_0N\|\rho\|_X$ for some uniform constant $C_0$.
\end{lemma}
\begin{proof}
Continuity of $w$ follows from the regularity  given in Lemma \ref{lem: w basic}. To get the support estimate, note that 
\begin{align*}
\left|\frac1{|\cdot|}*\rho(x)\right|&\le \|\rho\|_X\int_{\real^3}\frac{1}{|y|\langle x-y\rangle ^4}~dy\\
&\le \|\rho\|_X\left(\int_{|y|>\frac{|x|}2}\frac{2}{|x|\langle x-y\rangle^4}~dy+\int_{|y|<\frac{|x|}{2}}\frac{1}{|y|\langle x/2\rangle ^4}~dy\right)\\
&\le \frac{C_0\|\rho\|_X}{|x|}.
\end{align*}
So if $|x|> C_0N\|\rho\|_X$, $\left|\frac1{|\cdot|}*\rho(x)\right|<\frac1N$. By the definition of $\O_N$, we have $\alpha<-\frac1N$. The result now follows from the fact that $w(\kappa,r,u)=0$ when $u<0$, as stated in Lemma \ref{lem: w basic}.
\end{proof}

\begin{lemma}  
 The basic operator $\F:\mathcal O_\infty\to X\times\mathbb R$ 
 is $C^1$ Fr\'echet differentiable.  Its  Fr\'echet derivative is 
\begin{equation}
\frac{\partial\mathcal F}{\partial(\rho,\alpha,\kappa)}(\delta\rho,\delta\alpha,\delta\kappa) =\left(\delta\rho-\mathcal L(\delta\rho,\delta\alpha,\delta\kappa),\int_{\mathbb R^3} \delta\rho(x)~dx\right),
\end{equation}
with
\begin{align}   \label{F derivative}
\mathcal L(\delta\rho,\delta\alpha,\delta\kappa) = &~\partial_\kappa w\left(\kappa, r, \frac1{|\cdot|}*\rho(x) + \alpha\right) \cdot \delta\kappa\notag\\
&~+ \partial_u w\left(\kappa, r, \frac1{|\cdot|}*\rho(x) + \alpha\right)\cdot\left[\frac1{|\cdot|}*\delta\rho(x) + \delta\alpha\right].
\end{align}
\end{lemma}

\begin{proof}
We may fix any $N$ and restrict $\F$ to the subset  $\O_N$.   
Its first component $\F_1$ obviously maps $\O_N$ into $X$ 
because $w(\kappa,r,\frac1{|\cdot|}*\rho+\alpha)$ is continuous and compactly supported on a ball of radius $C_0N\|\rho\|_X$ by Lemma \ref{lem: comp supp}. For $\rho$ in a bounded set in $X$, $w(\kappa,r,\frac1{|\cdot|}*\rho+\alpha)$ is supported in a ball of fixed radius $R$. Since $\|f\|_X\le \langle R\rangle^4 \|f\|_{C(\overline{B_R})}$ for $f$ is supported in $\overline{B_R}$, it suffices to show Fr\'echet differentiability of the mapping into $C(\overline{B_R})\times \real$.  This is a simple consequence of the $C^1$ regularity of $w$.
\end{proof}

Our first goal is to find solutions to $\F(\rho,\alpha,\kappa)=0$ that are close to 
$(\rho_0,\alpha_0,0)$. 

\begin{lemma}\label{lem: iso}  
Given $a_0\in\A$, let $(\rho_0,\alpha_0)=(\rho(\cdot,a_0),\alpha(a_0))$ 
as in Lemma \ref{lem: radial uniqueness}.   
Then the operator $\L_0 := 
\frac{\partial\mathcal F}{\partial(\rho,\alpha)}(\rho_0,\alpha_0,0) : 
X \times\mathbb R \to X \times\mathbb R$ 
is an isomorphism provided that  \eqref{cond: M 1} is satisfied. 
\end{lemma}

\begin{proof}
It is easily seen that $\L_0$ is a compact perturbation of the identity, and it is  thus a Fredholm operator of index $0$. 
By Lemma \ref{lem: radial uniqueness}, $M_0$ in \eqref{eq: F2} is equal to $M(a_0)$. 
Thus it suffices to show that $\L_0$ is injective. By definition, $G(u) = w(0,0,u)$. 
By \eqref{F derivative}, 
\begin{equation}
\delta\rho(x)  - G'\left(\frac1{|\cdot|}*\rho_0(x) + \alpha_0\right)\cdot\left[\frac1{|\cdot|}*\delta\rho(x) + \delta\alpha\right]=0,
\end{equation}
\begin{equation}
\int_{\mathbb R^3}\delta\rho(x)~dx=0.  
\end{equation}
We have to show that $\delta\rho=0$, $\delta\alpha=0$. This can be proven in exactly 
the same way as Lemma 4.3 of \cite{strauss2019rapidly}, with $h^{-1}$ replaced by $G$.  
Indeed, let $u_0= \frac1{|\cdot|}*\rho+\alpha$ and $w = \frac1{|\cdot|}*\delta\rho+\delta\alpha$.  
We show that $\Delta w= -4\pi \frac{\rho_0'}{u_0'}w$ in the ball $B_0$ of radius $R(u_0)$, 
and $w=0$ outside $B_0$.  In \cite{strauss2019rapidly} it is proven first 
that $w$ is radial and then by a delicate argument that $w\equiv0$.  
Thus $\delta\rho\equiv0$ and $\delta\alpha=0$.  
\end{proof}
Lemma \ref{lem: iso} and the implicit function theorem imply the existence of a local curve of solutions near $(\rho_0,\alpha_0,0)$. 
However, to continue this solution curve globally, we use a global implicit function theorem that relies on the Leray-Schauder degree. To that end, we need
\begin{lemma}
The nonlinear operator $\N: (\rho(x),\alpha,\kappa)\mapsto w\left(\kappa,r,\frac1{|\cdot|}*\rho(x)+\alpha\right)$ is a compact map from $\O_\infty$ to $X$.
\end{lemma}
\begin{proof}
Call this nonlinear operator $\mathcal N$. By definition (see \cite{nirenberg1974topics}), 
we only have to show that $\mathcal N(\K)$ is precompact for every closed bounded 
set $\K\subset \O_\infty$.    Every such $\K$ is contained in some $\O_N$ 
for a finite $N$.  So by Lemma \ref{lem: comp supp}, the functions in $\mathcal N(\K)$ 
are supported on some fixed ball $B_R$. 
One can again dominate the $X$ norm by the $C(\overline{B_R})$ norm, 
and compactness is now obvious by Ascoli-Arzela.  
\end{proof}

We now state the global implicit function theorem to be applied to this problem.

\begin{theorem}   \label{thm: GIFT}
Let $X$ be a Banach space and let $U$ be an open subset of $X\times\real$.  
Let $F:U\to X$ be a Fr\'echet  $C^1$ mapping.  
Let $(\xi_0, \kappa_0)\in U$ such that $F(\xi_0,\kappa_0)=0$.  
Assume that the linear operator $\frac{\pa F}{\pa\xi}(\xi_0,\kappa_0)$ is an isomorphism on $X$.  
Assume that the mapping $(\xi.\kappa) \mapsto F(\xi,\kappa)-\xi$ is compact from $U$ to $X$.  
Let $\S$ be the closure in $X\times\real$ of the solution set $\{(\xi,\kappa)\  |\ F(\xi,\kappa)=0\}$. 
Let $\K$ be the connected component of $\S$ to which $(\xi_0,\kappa_0)$ belongs.   
Then one of the following three alternatives is valid. 
\begin{enumerate}[(i)]
\item $\K$ is unbounded in $X\times\real$.

\item $\K\backslash \{(\xi_0,\kappa_0)\}$ is connected.

\item $\K \cap \pa U \ne \emptyset$.  
\end{enumerate}
\end{theorem}
\begin{proof}
This is a standard theorem basically due to Rabinowitz. Theorem 3.2 in \cite{rabinowitz1971some}
in the case that $U=X\times\real$ and under some extra structural assumption.  
A more general version also appears in Theorem II.6.1 of \cite{kielhofer2006bifurcation};   
its proof is easy to generalize to permit a general open set $U$.  
The case of a general open set $U$ also appears explicitly in \cite{alexander1976implicit}.  
\end{proof}

We apply Theorem \ref{thm: GIFT} to $\F$ on $\O_\infty$ to obtain a preliminary form of Theorem \ref{main thm}:
\begin{lemma}\label{lem: pre main}
Assuming the mass condition \eqref{cond: M 2}, there is a connected set $\K$ of solutions to $\F(\rho,\alpha,\kappa)=0$ which contains $(\rho_0,\alpha_0,0)$ 
such that at least one of the following three alternatives is true:
\begin{enumerate}[(a)]
\item $\sup_{(\rho,\alpha,\kappa)\in\K} \|\rho\|_\infty=\infty$.
\item $\sup_{(\rho,\alpha,\kappa)\in\K} \text{diam}(\text{supp }\rho)=\infty$.
\item $\sup_{(\rho,\alpha,\kappa)\in\K} |\kappa| = \infty$, and there are two negative constants $a<b<0$ such that $a\le \alpha\le b$ for all $(\rho,\alpha,\kappa)\in \K$.
\end{enumerate}
\end{lemma}

\begin{proof}
We apply Theorem \ref{thm: GIFT} to $F=\F$, $X=Z=X\times \real$, $U=\O_N$, and $\xi=(\rho,\alpha)$, with starting point $(\xi_0,\kappa_0)=(\rho_0,\alpha_0,0)$.   
Thus there exists a solution set $\K_N$ such that at least one of the three alternatives in Theorem \ref{thm: GIFT} holds. 
We  claim that alternative (ii) in Theorem \ref{thm: GIFT} cannot happen.  
Indeed, let $\K_N\setminus \{(\rho_0,\alpha_0,0)\}$ be connected.  
(This means that $\K_N$ contains a ``loop".)   
Since the projection onto the $\kappa$-axis is continuous and $\K_N$ obviously must contain  solutions with $\kappa$ positive and negative, 
$\K_N\setminus \{(\rho_0,\alpha_0,0)\}$ must also contain a solution of the form $(\rho_1,\alpha_1,0)$, where $\rho_1$ is compactly supported by Lemma \ref{lem: comp supp} and $\F_1(\rho_1,\alpha_1,0)=0$. Lemma \ref{lem: non-rot} implies that $\rho_1$ must be one of the radial solutions supported on a compact ball. Since $\rho_1$ has the same total mass as $\rho_0$ by $\F_2(\rho_1,\alpha_1,0)=0$, the mass condition \eqref{cond: M 2} implies that $\rho_1(0)=\rho_0(0)$. Lemma \ref{lem: radial uniqueness} implies that $(\rho_1,\alpha_1)=(\rho_0,\alpha_0)$. This is a contradiction.
Thus, only cases (i) or (iii) in Theorem \ref{thm: GIFT} can happen. 

In other words, either
$$
\sup_{(\rho,\alpha,\kappa)\in\K_N}(\|\rho\|_X+|\kappa|+|\alpha|)=\infty,
$$
or 
$$
\sup_{(\rho,\alpha,\kappa)\in\K_N} \alpha = -\frac1N.
$$
Now the set $\K=\cup_{N=1}^\infty \K_N$ is connected because each $\K_N$'s is connected and they are nested.   For $\K$ we obviously have either 
\eqn\label{eq: sup blow up}
\sup_{(\rho,\alpha,\kappa)\in\K}(\|\rho\|_X+|\kappa|+|\alpha|)=\infty,
\eeqn
or 
\eqn\label{eq: sup alpha}
\sup_{(\rho,\alpha,\kappa)\in\K} \alpha = 0.
\eeqn

Suppose now that neither (a) nor (b) in the statement of the lemma occurs. Then there exists a constant $H>0$ such that for all $(\rho,\alpha,\kappa)\in\K$, 
both $\|\rho\|_\infty\le H$ and  supp $\rho$ is contained in the ball centered at the origin with radius $H$. We claim that \eqref{eq: sup alpha} cannot occur. In fact, for any $\bar x$ with $|\bar x|=2H$, we have the lower bound 
\eqn 
\frac1{|\cdot|}*\rho(\bar x)\ge \frac{1}{|\bar x|+H}\int_{|y|\le C}\rho(y)~dy=\frac{M_0}{3H}. 
\eeqn 
So if there were to exist a point $(\rho,\alpha,\kappa)\in \K$ with $\alpha>-\frac{M_0}{3H}$, we would have $\frac1{|\cdot|}*\rho(\bar x) + \alpha >0$ and 
$$
\rho(\bar x) = w\left(\kappa, r(\bar x), \frac1{|\cdot|}*\rho(\bar x) + \alpha\right)>0.
$$
This contradicts the assumption that $\rho$ is supported in $|x|<H$. 

Thus \eqref{eq: sup blow up} must occur, and $b:=\sup_{(\rho,\alpha,\kappa)\in\K} \alpha<0$. On the other hand, for $|x|\le H$, we have the upper bound
\eqn
0\le \frac1{|\cdot|}*\rho(x)\le \int_{|y|\le H}\frac H{|x-y|}~dy\le H\int_{|y|\le 2H}\frac1{|y|}~dy=8\pi H^3.
\eeqn
If $\alpha< -8\pi H^3$, then $\rho(x)=0$ for all $|x|\le H$, which implies that $\rho$ is identically zero. This contradicts the mass constraint $\F_2(\rho,\alpha,\kappa)=0$. 
Thus we can define $a=-8\pi H^3$ and we have $a\le\alpha\le b$ for any $(\rho,\alpha,\kappa)\in \K$. We also obviously have $\|\rho\|_X\le \langle H\rangle^4  H$. 
So the only way for \eqref{eq: sup blow up} to occur is to have
$$
\sup_{(\rho,\alpha,\kappa)\in\K} |\kappa| = \infty.
$$
This means that (c) occurs.  
\end{proof}

Any solution $(\rho,\alpha,\kappa)\in \K$ satisfies $\rho\in C^1$ because 
$\rho$ is continuous, $\frac1{|\cdot|}*\rho \in C^1$ 
and $\rho = w(\kappa, r, \frac1{|\cdot|}*\rho + \alpha)$.  

By Lemma 4.5, we are left with one final possibility to be eliminated in order to prove Theorem \ref{main thm}. That possibility is: $\|\rho\|_\infty$ and the support of $\rho$ are both uniformly bounded, $\alpha$ is bounded between two negative numbers, and $\kappa$ is unbounded. It turns out that a much more delicate argument is needed to eliminate this last possibility.  It is provided in the next two sections.  
\section    {Eliminate unbounded $\kappa$ for a special example}\label{sec: special}

To better explain the essential ideas of the argument, we prove Theorem \ref{main thm} in this section for the special example of $\phi$ given as 
\eqn
\phi(E,L) = (-E)_+^{\frac12} (C_1+C_2L^2)
\eeqn
for suitable constants $C_1,C_2>0$.
The proof for the general case will be given in the next section.  By a direct calculation, we get for $u>0$:
\begin{align}
w(\kappa, r,u)& = 2\pi\int_{-u}^0 \int_{-\sqrt{2(E+u)}}^{\sqrt{2(E+u)}} (-E)_+^{\frac12} [C_1+C_2(\kappa r s)^2]~ds~dE\notag\\
&= 4\pi C_1u^2\int_{-1}^0 [2(E+1)(-E)]^{\frac12}~dE+\frac{4\pi}3 C_2\kappa ^2 r^2 u^3\int_{-1}^0[8(E+1)^3(-E)]^{\frac12}~dE\notag\\
&=u^2+\kappa^2 r^2 u^3.
\end{align}
Here we have chosen $C_1$ and $C_2$ to make the coefficients equal to 1.
Thus we may write $\F_1(\rho,\alpha,\kappa)=0$ as
\begin{align}   \label{eq: model case}
\rho &= \kappa^2 r^2\left(\frac1{|\cdot|}*\rho+\alpha\right)_+^3+\left(\frac1{|\cdot|}*\rho+\alpha\right)_+^2\notag\\
&=\kappa^2 r^2 u_+^3+u_+^2.    
\end{align}
Under the above setup, 
we want to prove Theorem \ref{main thm} as a consequence of the following lemmas. Theorem \ref{main thm} means that either case (a) or case (b) of Lemma \ref{lem: pre main} must occur. If neither case (a) nor case (b) occurs, we get from Lemma \ref{lem: pre main} that case (c) must occur, and there is a constant $R>0$ such that for all $(\rho,\alpha,\kappa)\in \K$, we have $\|\rho\|_\infty<R$, and supp $\rho$ is contained in $B_R$, the ball of radius $R$ centered at the origin.
%
%
%
%
Denote  $u = \frac1{|\cdot|}*\rho+\alpha$.   From \eqref{eq: model case} we see that 
$u>0$ in the region $\Omega =\{\rho>0\}$ occupied by the star,  
while $u<0$ outside $\overline{\Omega}$.  
First we show that $u$ must be small on $\Omega$ if $\kappa$ is large. In fact, we have
\begin{lemma}\label{lem: kappa -2/5}
There exists a $C>0$ such that  $\|u\|_{\infty;\Omega}\le C\kappa^{-2/5}$ for $|\kappa|>1$.
\end{lemma}
\begin{proof}
$u>0$ on $\Omega$ so $u=u_+$ there.  
Because of the product $\kappa^2 r^2$ in \eqref{eq: model case}, we have to treat small $r$ separately.  
Let $r_0>0$ be a radius to be chosen sufficiently small later on. 
For $x$ such that $r(x)\ge r_0$, we use \eqref{eq: model case} to get
\eqn   \label{first bound on u_+}
u_+(x)\le \sqrt[3]{\frac{\|\rho\|_\infty}{\kappa^2 r_0^2}}.
\eeqn
For $x$ such that $r(x)\le r_0$, we use the simple potential estimate $\|Du\|_\infty\le C\|\rho\|_\infty$ and integrate it radially from $r(x)$ to $r_0$.   
Thereby we obtain 
\eqn  \label{potential estimate}
u_+(x)\le \sqrt[3]{\frac{\|\rho\|_\infty}{\kappa^2 r_0^2}}+C\|\rho\|_\infty r_0 
\eeqn
for all $x\in\real^3$.  
We now choose $r_0=r_0(\kappa)$ so that the two terms above are comparable. As a result we get 
\eqn \label{bound on u_+}
\|u\|_{\infty;\Omega}=\|u_+\|_{\infty;\real^3} \le C\kappa^{-2/5} 
\eeqn
with $C$ independent of $\kappa$.  
\end{proof}

Note that we cannot obtain a smallness estimate for $u$ outside $\Omega$, because \eqref{eq: model case} loses control of $u$ when $u<0$. While $u$ is close to $\alpha$ near infinity, $\alpha$ is bounded but not necessarily small. Our second step is to control $D^2u$ by a low power of $\kappa$.  
\begin{lemma}
Let $B_R$ be the common ball containing the fluid domain $\Omega$ given above. For every $\tau\in(0,1)$, there exists a constant $C_\tau>0$ such that
\eqn\label{bound on D^2u}
\|D^2 u\|_{\infty;B_R}\le C_\tau\kappa^{2\tau}.
\eeqn
\end{lemma}
\begin{proof}
Indeed, we differentiate \eqref{eq: model case} to get
\eqn
D\rho = \kappa^2 Dr^2 u_+^3+\kappa^2 r^23u_+^2 Du + 2u_+ Du, 
\eeqn
so that
\eqn
\|D\rho\|_{\infty;B_R} \le C\kappa^2(1+\|Du\|_{\infty;B_R})\le C\kappa^2(1+\|\rho\|_{\infty;\real^3})\le C\kappa^2.
\eeqn
For any $\tau \in(0,1)$  we have the interpolation inequality
\eqn   \label{interp}
\|\rho\|_{C^{0,\tau}(B_R)}\le C_{\tau}\|\rho\|_{\infty;B_R}^{1-\tau}\|D\rho\|_{\infty;B_R}^\tau.
\eeqn
From the previous two estimates we get 
\eqn
\|\rho\|_{C^{0,\tau}(B_R)}\le C_{\tau}\kappa^{2\tau}.
\eeqn
Now by standard potential estimates we get
\eqn 
\|D^2 u\|_{\infty;B_R}  \le   C_\tau\|\rho\|_{C^{0,\tau}(B_R)}  
\le   C_{\tau}\kappa^{2\tau}.
\eeqn
\end{proof}
Our third step is to use an elementary Gagliardo-Nirenberg (GN) type of inequality in the fluid domain $\Omega$ in order to obtain smallness of $Du$ there. This kind of inequality is well-known, but the particular form we are using here is more difficult to find. We include the proof here for completeness.  
\begin{lemma}
Let $\Omega$ be a bounded open set in $\real^n$.  
Let $u\in C(\overline{\Omega})\cap C^2(\Omega)$ such that $u=0$ on $\pa\Omega$.  Then there exists a constant $C$ independent of $u$ and $\Omega$ such that
\eqn\label{eq: GN}
\|Du\|_{\infty;\Omega}\le C\|u\|_{\infty;\Omega}^{1/2}\ \|D^2 u\|_{\infty;\Omega}^{1/2}.   \eeqn
\end{lemma}

\begin{proof}   
The key point here is that the constant $C$ is uniform and independent of the domain $\Omega$, {\it no matter how rough its boundary might be}. The idea is to integrate the second derivative from a nearby point, and thereby estimate $Du$ at a nearby point by a difference quotient. It is crucial that $u$ vanishes on $\partial \Omega$ in order to handle an extremely close difference quotient.

Without loss of generality, we may assume $\|D^2u\|_{\infty;\Omega}\ne 0$, since if it does, then $u$ is a linear function on each connected component of $\Omega$, hence must be zero identically by the boundary condition.  
Define $d^2 = \|u\|_{\infty;\Omega}/\|D^2u\|_{\infty;\Omega}$.  Take any point $x\in\Omega$ 
and any direction $e_i$.  
In case the line segment $(x,x+de_i)\subset\Omega$, 
then there is a value $\xi\in (0,d)$ so that $u(x+de_i)-u(x) = d\cdot \pa_i u(x+\xi e_i)$.  
Hence $|\pa_i u(x+\xi e_i)| \le (2/d)\|u\|_{\infty;\Omega}$.  
Also, $|\pa_iu(x)-\pa_iu(x+\xi e_i)| = |\int _\xi^0 \pa_i^2 u(x+te_i) dt| \le d\|D^2u\|_{\infty;\Omega}$.  
Combining the last two inequalities, we have 
 $|\pa_iu(x)| \le \tfrac2d \|u\|_{\infty;\Omega} + d\|D^2u\|_{\infty;\Omega}= 3\sqrt{\|u\|_{\infty;\Omega} \|D^2u\|_{\infty;\Omega}}$.  
 The last equality is due to the definition of $d$.  
The same argument is valid in case the line segment $(x-de_i,x)\subset\Omega$.

The other possibility is that $x$ and $i$ are such that there is no such line segment.  
In that case consider the line segment $[x-ae_i,x+be_i]$ with both $a\le d$ and $b\le d$  
such that  both endpoints lie on $\pa\Omega$.
Since $u$ vanishes at both endpoints, there is an intermediate value $\xi$ such that $\pa_i u(x+\xi e_i)=0$.  
Thus $|\pa_iu(x)| = |\int _\xi^0 \pa_i^2 u(x+te_i) dt| \le d\|D^2u\|_{\infty;\Omega} 
= \sqrt{\|u\|_{\infty;\Omega} \|D^2u\|_{\infty;\Omega}}$.   Summing on $i\in\{1,2,3\}$, we obtain \eqref{eq: GN}.  
\end{proof}

We can now combine the previous inequalities to get a smallness estimate for $Du$ on $\Omega$.  
Into the GN inequality \eqref{eq: GN} we substitute the bound \eqref{bound on u_+} on $u$ 
and the bound \eqref{bound on D^2u}  on $D^2u$ to obtain 
\eqn
\|Du\|_{\infty,\Omega}\le C_\tau|\kappa|^{-1/5}|\kappa|^{\tau} 
\le C|\kappa|^{-\frac1{10}}      \eeqn
for large $\kappa$, where we have chosen $\tau=1/{10}$. 

If $\Omega$ is sufficiently smooth, we can directly integrate $Du$ on $\partial\Omega$ and obtain a smallness estimate for $\|\rho\|_1$. However, since $\partial\Omega$ could potentially be very rough or have large measure, we need to take a final step to extend the smallness estimate of $Du$ to the whole space. Fortunately this is easy to do, due the harmonicity of $Du$ outside $\Omega$. In fact $\Delta u = -4\pi\rho=0$, so $u$ and as consequence $Du$ are harmonic outside $\Omega$.  
Also $\lim_{x\to\infty} Du=0$. By the maximum principle, 
\eqn
\|Du\|_{\infty, \real^3}  \le  C|\kappa|^{-\frac1{10}}.
\eeqn
Now we integrate $\rho =-\tfrac1{4\pi}\Delta u$ on the ball $B_R$ to get
\eqn \label{M vs kappa}
M_0 = \int_{B_R} \rho\ dx = -\frac1{4\pi} \int_{\partial B_R} \nabla u \cdot n ~d\sigma\le C|\kappa|^{-\frac1{10}}.
\eeqn
Letting $|\kappa|\to\infty$, we find $M_0=0$, which is a contradiction. This finishes the proof of  Theorem \ref{main thm}.


\section{Eliminate unbounded $\kappa$ for general $\phi$}

We impose general assumptions on $\phi(E,L)$ that are given in Section 2 
and then repeat the preceding proof.
 We now turn to the general case. Recall that $\phi$ satisfies \eqref{phi: lower bound} and \eqref{phi: upper bound}. We will proceed as in Section \ref{sec: special}. The next lemma generalizes \eqref{first bound on u_+}.  
\begin{lemma}\label{lem: w lower bound}
For every $H>0$, there exist $C>0$, $K>0$, such that for 
any triple $(\kappa,r,u)$  satisfying 
\eqn  \label{eq: three est}
|w(\kappa,r,u)|\le H  , \quad  u>0  , \quad  |\kappa r |>K,   \eeqn
we have 
\eqn
u\le C|\kappa r|^{-\frac{2\Gamma}{3+\Gamma+2\delta}}.   \eeqn
Here $\Gamma$ and $\delta$ are the numbers given in \eqref{phi: lower bound}.
\end{lemma}

\begin{proof}
By \eqref{phi: lower bound}, there exist $C>0$, $A>0$, $\Gamma>0$, $\delta>0$ such that for $-A<E<0$, $|L|>1/A$ we have the lower bound  
\eqn\label{eq: phi lower bound}
\phi(E,L)\ge C|E|^{\delta}|L|^\Gamma,   \eeqn
For simplicity of notation in this proof, we take $A=1$.  
Note that in the following integrations, we can only use \eqref{eq: phi lower bound} when $-1<E<0$, and $|\kappa r s|>1$, or $|s|>\frac{1}{|\kappa r|}$. 

We claim that $u\le 1$ if \eqref{eq: three est} holds. Indeed, in case $u>1$, we have
\begin{align*}
w(\kappa,r,u)&\ge C\int_{-\frac12}^0\int_{-\sqrt{2(E+1)}}^{\sqrt{2(E+1)}}\phi(E,\kappa r s)~ds~dE\\
&\ge C\int_{-\frac12}^0\int_0^1 \phi(E,\kappa r s)~ds~dE.
\end{align*}
We pick $K>2$ so large that $\frac1{|\kappa r|}<\frac 1K<\frac12$. 
It follows that
\begin{align*}
w(\kappa,r,u)&\ge C\int_{-\frac12}^0\int_{\frac1{|\kappa r|}}^1\phi(E,\kappa r s)~ds~dE\\
&\ge C\int_{-\frac12}^0\int_{\frac1{|\kappa r|}}^1~|E|^\delta |\kappa r |^\Gamma s^\Gamma ds~dE\\
&\ge C|\kappa r|^\Gamma \int_{-\frac12}^0\int_{\frac12}^1~|E|^\delta s^\Gamma ds~dE\\
&=C|\kappa r|^\Gamma.
\end{align*}
We also pick $K>\left(\frac H C\right)^{\frac1\Gamma}$, so that the 
assumption $|\kappa r|>K$ and the preceding inequality imply $w(\kappa, r, u)>H$. 
Since this contradicts the first inequality in \eqref{eq: three est}, it follows that $0<u<1$.  

In case $\frac{1}{|\kappa r|}>\frac{\sqrt u}{2}$, we have 
\eqn
u<4|\kappa r|^{-2}<C|\kappa r|^{-\frac{2\Gamma}{3+\Gamma+2\delta}},     \eeqn
as desired, 
because $\frac{2\Gamma}{3+\Gamma+2\delta}<2$.   
On the other hand, in case $\frac1{|\kappa r|}\le \frac{\sqrt u}{2}$,  we have
\begin{align*}
w(\kappa, r, u)&\ge C\int_{-\frac u2 }^0\int_{\frac1{|\kappa r|}}^{\sqrt{2(E+u)}}\phi(E,\kappa r s)~ds~dE\\
&\ge C\int_{-\frac u 2}^0\int_{\frac{\sqrt u}2}^{\sqrt u}|E|^\delta |\kappa r|^\Gamma s^\Gamma ~ds~dE\\
&\ge C|\kappa r|^\Gamma u^{\frac{\Gamma+1}{2}+\delta+1}.
\end{align*}
Thus the hypothesis $w(\kappa, r, u)\le H$ implies
\eqn
u\le \left(\frac{H}{C|\kappa r|^\Gamma}\right)^{\frac{2}{3+\Gamma+2\delta}}\le C|\kappa r|^{-\frac{2\Gamma}{3+\Gamma+2\delta}}.
\eeqn
\end{proof}

Since the inner integration of the definition of $w$ always traverses the symmetric interval $(-\sqrt{2(E+u)},\sqrt{2(E+u)})$, the preceding proof can be modified to work by assuming the limit inequality \eqref{phi: lower bound} only as $L\to\infty$ or only as $L\to -\infty$.

The next lemma will be used to replace the simple powers in \eqref{eq: model case}. 
\begin{lemma}\label{lem: w upper bound}
For every $B>0$, there exists $C>0$ such that for $|r|\le B$, $|u|\le B$ and $|\kappa|>1$, 
we have the upper bounds 
\eqn
|\partial_r w(\kappa,r,u)|\le C |\kappa|^{1+\Lambda}.\eeqn 
and 
\eqn
|\partial_u w(\kappa,r,u)|\le C |\kappa|^\Lambda.\eeqn
Here $\Lambda$ is the number given in \eqref{phi: upper bound}.
\end{lemma}

\begin{proof}
The derivatives of $w$ are 
\eqn
\partial_rw(\kappa,r,u)=2\pi\kappa\int_{-u}^0\int_{-\sqrt{2(E+u)}}^{\sqrt{2(E+u)}}\partial_L\phi(E,\kappa rs)s~ds~dE.
\eeqn
\eqn
\partial_uw(\kappa,r,u)=\sqrt2\pi\int_{-u}^0\left(\phi(E,\kappa r \sqrt{2(E+u)})+\phi(E,-\kappa r \sqrt{2(E+u)})\right)\frac{dE}{\sqrt{E+u}}.
\eeqn
By assumption \eqref{phi: upper bound} and $u>0$, the $r$-derivative has the bound 
\begin{align*}
|\partial_r w (\kappa,r,u)|&\le C|\kappa |\int_{-u}^0\int_{-\sqrt{2(E+u)}}^{\sqrt{2(E+u)}}  
|E|^{-\mu} (1+|\kappa r|^\Lambda |s|^\Lambda)|s|~ds~dE\\
&\le C|\kappa|\int_{-u}^0|E|^{-\mu} 
\left[(E+u)+|\kappa r|^\Lambda (E+u)^{(2+\Lambda)/2}\right]~dE\\
&\le C|\kappa|u^{2-\mu}\int_{-1}^0|F|^{-\mu}(1+F)~dF+C|\kappa |^{1+\Lambda} u^{{\frac{4+\Lambda-2\mu}{2}}}\int_{-1}^0|F|^{-\mu}(1+F)^{\frac{2+\Lambda}{2}}~dF\\
&\le C(|\kappa|+|\kappa|^{1+\Lambda})
\end{align*}
for $|u|<B$, $r<B$, because $\mu<\frac12$.  
By assumption \eqref{phi: upper bound}, the $u$-derivative has the bound 
\begin{align*}
|\partial_u w(\kappa,r,u)|&\le C\int_{-u}^0|E|^{-\mu}(1+|\kappa r|^\Lambda(E+u)^{\Lambda/2})(E+u)^{-1/2}~dE\\
&\le C u^{\frac12-\mu}\int_{-1}^0 |F|^{-\mu}(1+F)^{-1/2}~dF+C|\kappa|^\Lambda u^{\frac{\Lambda+1-2\mu}{2}}\int_{-1}^0|F|^{-\mu}(1+F)^{(\Lambda-1)/2}~dF\\
&\le C(1+|\kappa|^\Lambda)
\end{align*}
for $|u|<B$, $r<B$, because $\mu<\frac12$.
\end{proof}

Using the above two lemmas, we can complete the proof as in the special case of 
Section \ref{sec: special}. 
Indeed, we have $\rho=w(\kappa,r,u)$, $u=\frac1{|\cdot|}*\rho+\alpha$, and $\int \rho~dx=M_0$. Let there be a constant $R>0$ such that $\rho$ is supported in $B_R$, and $\|\rho\|_\infty\le R$, $-R<\alpha<-\frac1R$. We want to show that $\kappa$ must be bounded. 
On the contrary, suppose $\sup_{(\rho,\alpha,\kappa)\in\K}|\kappa|=\infty$. 

Let $\Omega=\{u>0\}$ as before.  First, we use the standard potential estimates to get $\|Du\|_{\infty}\le C\|\rho\|_\infty\le C$.    Define 
\eqn
r_0 = r_0(\kappa) = |\kappa|^{-\frac{2\Gamma }{3+3\Gamma+2\delta}}. 
\eeqn
If $r(x)>r_0$, then $|\kappa r|>|\kappa|^{\frac{3+\Gamma+2\delta}{3+3\Gamma+2\delta}}$. Apply Lemma \ref{lem: w lower bound} with $H=\|\rho\|_\infty$. Since $w(\kappa,r,u)=\rho<H$, there are constants $C$ and $K$ such that if $|\kappa|^{\frac{3+\Gamma+2\delta}{3+3\Gamma+2\delta}}>K$,
\eqn    \label{2bound on u}
u(x)\le C|\kappa r|^{-\frac{2\Gamma}{3+\Gamma+2\delta}}\le C|\kappa|^{-\frac{2\Gamma}{3+3\Gamma+2\delta}}.
\eeqn
If $r(x)\le r_0$, we integrate from the cylinder of radius $r_0$ and use the bound 
$\|Du\|_\infty \le C$ to get
\eqn
\|u\|_{\infty;\Omega}\le C|\kappa|^{-\frac{2\Gamma}{3+3\Gamma+2\delta}}+Cr_0\le C|\kappa|^{-\frac{2\Gamma}{3+3\Gamma+2\delta}}.
\eeqn 
This generalizes Lemma \ref{lem: kappa -2/5}.  

Secondly, we differentiate the equation $\rho=w(\kappa,r,u)$ to get
\eqn
D\rho = \partial_rw(\kappa,r,u)Dr + \partial_u w(\kappa,r,u)Du.
\eeqn
Using Lemma \ref{lem: w upper bound} and the bound $\|Du\|_\infty\le C$ we get
\eqn
\|D\rho\|_{\infty;B_R}\le C|\kappa|^{1+\Lambda}.
\eeqn
By interpolation \eqref{interp}, we obtain for $\tau\in (0,1)$
\eqn
\|\rho\|_{C^{0,\tau}(B_R)} \le C_{\tau}\|\rho\|_\infty^{1-\tau}\|D\rho\|_{\infty;B_R}^\tau
\eeqn
and the bound on $\|\rho\|_\infty$, 
\eqn
\|\rho\|_{C^{0,\tau}(B_R)} \le C_{\tau}|\kappa|^{\tau(1+\Lambda)}.
\eeqn
By standard potential estimates on $B_R$, 
\eqn   \label{2bound on D^2u}
\|D^2 u\|_{\infty;B_R}\le C_\tau|\kappa|^{\tau(1+\Lambda)}.
\eeqn
Thirdly, using the GN inequality \eqref{eq: GN} on $\Omega$ exactly as before
together with the bounds \eqref{2bound on u}  and  \eqref{2bound on D^2u}, 
we obtain for $|\kappa|$ sufficiently large
\eqn
\|Du\|_{\infty;\Omega}\le C_\tau |\kappa|^{\frac{\tau(1+\Lambda)}2-\frac{\Gamma}{3+3\Gamma+2\delta}}.
\eeqn
Finally, $Du$ is harmonic outside $\Omega$ and vanishes at infinity, 
so that
\eqn
\|Du\|_{\infty;\real^3}\le C_\tau |\kappa|^{\frac{\tau(1+\Lambda)}2-\frac{\Gamma}{3+3\Gamma+2\delta}}.
\eeqn
We  choose $\tau$ small enough so that the last exponent is negative. 
By integrating $Du$ over $\partial B_R$ and letting $|\kappa|\to\infty$, 
we deduce that the mass $M_0$ vanishes, which is a contradiction. This finishes the proof of Theorem \ref{main thm}.
\qed 

\section 
{Eliminate unbounded $\rho$ for some $\phi$}
The purpose of this short section is to prove Theorem \ref{main thm2}. 
For this theorem we do not require hypothesis \eqref{phi: lower bound}, but must  strengthen
\eqref{phi: upper bound} to require $\Lambda<4$. 
We first estimate the growth rate of $w$ respect to $u$.

\begin{lemma} \label{lem: w growth rate}
For every $B>0$, there exists a constant $C$ such that 
\eqn
|w(\kappa,r,u)|  \le  C (1+u_+^{1+\Lambda/2})    \eeqn 
for all $|\kappa|<B,r<B$. 
\end{lemma}
\begin{proof}
Recalling the definition of $w$ and using \eqref{phi: upper bound},  we have for $u>0$
\begin{align*}
|w(\kappa,r,u)|  \le  &C \int_{-u}^0  \int_{-\sqrt{2(E+u)}} ^{\sqrt{2(E+u)}} ~ 
|\phi(E,\kappa rs)| ~ ds~dE   \\
& \le C\int_{-u}^0  \int_{-\sqrt{2(E+u)}} ^{\sqrt{2(E+u)}} ~  |E|^{-\frac12} (1+|s|^\Lambda) ~ds~dE  \\
& \le C\int_{-u}^0  |E|^{-\frac12}  \left[  (E+u)^{\frac12} + (E+u)^{\frac{\Lambda+1}2 }\right]  dE  \\
& \le C\int_{-1}^0 (uF)^{-\frac12} \left[ u^{\frac12}(1+F)^{\frac12} 
+  u^{\frac{\Lambda+1}2}(1+F)^{\frac{\Lambda+1}2}\right]  u~dF   \\
& \le Cu + Cu^{1+\frac\Lambda 2}  .
\end{align*} 
Recall that $w=0$ for $u\le 0$.   
\end{proof}

{\it Proof of Theorem  \ref{main thm2}.}  
The proof is similar to that of Theorem 7.1 in \cite{strauss2019rapidly}. If the conclusion of Theorem \ref{main thm2} does not hold, then neither case (b) nor case (c) of Lemma \ref{lem: pre main} occurs. Thus case (a) must occur.   
A contradiction would follow from a uniform $L^\infty$ bound on $\rho$, assuming that the support of $\rho$ is contained in a fixed ball $B_R$, and that $|\kappa|$ is uniformly bounded.

Recalling that $\rho=w(\kappa,\frac1{|\cdot|} * \rho+\alpha)$ 
and writing $1+\frac\Lambda 2 = 3-\ep$ with $0<\ep<2$, we have by Lemma \ref{lem: w growth rate}
\eqn \label{eq: Lq improve}
\rho \le C\left[ 1 +  \left( \frac1{|\cdot|} * \rho+\alpha \right)_+^{3-\ep}\right]  
      \le C + C\left( \frac1{|\cdot|} * \rho \right)^{3-\ep},
\eeqn
since $\alpha<0$.
Now the total mass $M_0=\int\rho ~dx$ is fixed 
and we are assuming that $\rho$ has its support contained in $B_R$.  
By the endpoint Hardy-Littlewood-Sobolev inequality, $ \frac1{|\cdot|} * \rho$ is uniformly bounded in  $L^q(B_R)$ for every $1<q<3$. It follows from \eqref{eq: Lq improve} that $\rho$ is uniformly bounded in $L^{q/(3-\epsilon)}(B_R)$.
This is a slight improvement of the $L^1$ bound if we choose $q$ sufficiently close to 3.  
Repeated use of the Hardy-Littlewood-Sobolev inequality allows us to improve an $L^p$ bound to an $L^r$ bound on $\rho$, with $\frac1r=(3-\ep)\left(\frac1p-\frac23\right)$ for any $p\in (1,\frac32)$. It is not hard to see that a finite number of repetitions of such estimates 
leads to $\rho$ being bounded uniformly in $L^r$ for some $r>\frac32$, at which point we  conclude that $\frac1{|\cdot|}*\rho$ and consequently $\rho$ are uniformly bounded in $L^\infty(B_R)$.  \qed

\bibliographystyle{acm}
\bibliography{rotstarbiblio}

\begin{thebibliography}{10}

\bibitem{alexander1976implicit}
{\sc Alexander, J., and Yorke, J.~A.}
\newblock The implicit function theorem and the global methods of cohomology.
\newblock {\em Journal of Functional Analysis 21}, 3 (1976), 330--339.

\bibitem{ambrosetti1973dual}
{\sc Ambrosetti, A., and Rabinowitz, P.~H.}
\newblock Dual variational methods in critical point theory and applications.
\newblock {\em Journal of functional Analysis 14}, 4 (1973), 349--381.

\bibitem{andreasson2014rotating}
{\sc Andr{\'e}asson, H., Kunze, M., and Rein, G.}
\newblock Rotating, stationary, axially symmetric spacetimes with collisionless
  matter.
\newblock {\em Communications in Mathematical Physics 329}, 2 (2014), 787--808.

\bibitem{binney2011galactic}
{\sc Binney, J., and Tremaine, S.}
\newblock {\em Galactic dynamics}, vol.~13.
\newblock Princeton university press, 2011.

\bibitem{de1982priori}
{\sc De~Figueiredo, D.~G., Lions, P., and Nussbaum, R.}
\newblock A priori estimates and existence of positive solutions of semilinear
  elliptic equations.
\newblock In {\em Djairo G. de Figueiredo-Selected Papers}. Springer, 1982,
  pp.~133--155.

\bibitem{gidas1979symmetry}
{\sc Gidas, B., Ni, W.-M., and Nirenberg, L.}
\newblock Symmetry and related properties via the maximum principle.
\newblock {\em Communications in Mathematical Physics 68}, 3 (1979), 209--243.

\bibitem{guo2008unstable}
{\sc Guo, Y., and Lin, Z.}
\newblock Unstable and stable galaxy models.
\newblock {\em Communications in Mathematical Physics 279}, 3 (2008), 789--813.

\bibitem{jang2017slowly}
{\sc Jang, J., and Makino, T.}
\newblock On slowly rotating axisymmetric solutions of the {E}uler--{P}oisson
  equations.
\newblock {\em Archive for Rational Mechanics and Analysis 225}, 2 (2017),
  873--900.

\bibitem{kielhofer2006bifurcation}
{\sc Kielh{\"o}fer, H.}
\newblock {\em Bifurcation theory: An introduction with applications to PDEs},
  vol.~156.
\newblock Springer Science \& Business Media, 2006.

\bibitem{nirenberg1974topics}
{\sc Nirenberg, L.}
\newblock {\em Topics in nonlinear functional analysis}, vol.~6.
\newblock American Mathematical Soc., 1974.

\bibitem{rabinowitz1971some}
{\sc Rabinowitz, P.~H.}
\newblock Some global results for nonlinear eigenvalue problems.
\newblock {\em Journal of Functional Analysis 7}, 3 (1971), 487--513.

\bibitem{rein2007collisionless}
{\sc Rein, G.}
\newblock Collisionless kinetic equations from astrophysics--the
  {V}lasov-{P}oisson system.
\newblock {\em Handbook of differential equations: evolutionary equations 3\/}
  (2007), 383--476.

\bibitem{rein2003stable}
{\sc Rein, G., and Guo, Y.}
\newblock Stable models of elliptical galaxies.
\newblock {\em Monthly Notices of the Royal Astronomical Society 344}, 4
  (2003), 1296--1306.

\bibitem{rein2000compact}
{\sc Rein, G., and Rendall, A.~D.}
\newblock Compact support of spherically symmetric equilibria in
  non-relativistic and relativistic galactic dynamics.
\newblock In {\em Mathematical Proceedings of the Cambridge Philosophical
  Society\/} (2000), vol.~128, Cambridge University Press, pp.~363--380.

\bibitem{strauss2017steady}
{\sc Strauss, W.~A., and Wu, Y.}
\newblock Steady states of rotating stars and galaxies.
\newblock {\em SIAM Journal on Mathematical Analysis 49}, 6 (2017), 4865--4914.

\bibitem{strauss2019rapidly}
{\sc Strauss, W.~A., and Wu, Y.}
\newblock Rapidly rotating stars.
\newblock {\em Communications in Mathematical Physics 368}, 2 (2019).

\end{thebibliography}

\end{document}